\documentclass{amsart}
\newcommand{\nn}{\nonumber}
\newcommand{\ba}{\begin{eqnarray}}
\newcommand{\ea}{\end{eqnarray}}

\newcommand{\no}{\noindent}

\newcommand{\stone}{\left[ {m+1} \atop {j} \right]}
\newcommand{\sttwo}{\left[ {m+1} \atop {m} \right]}
\newcommand{\stthree}{\left[ {m+1} \atop {m+1-k} \right]}
\newcommand{\stfour}{\left[ {m} \atop {k} \right] }
\newcommand{\stfive}{\left[ {m} \atop {m-k} \right] }
\newcommand{\stsix}{\left[ {m} \atop {m-1} \right] }
\newcommand{\stseven}{\left[ {m} \atop {m-2} \right] }
\newcommand{\steight}{\left[ {m} \atop {m-3} \right] }
\newcommand{\stnine}{\left[ {m} \atop {m-4} \right]}
\newcommand{\stten}{\left[ {m+1} \atop {m-1} \right] }
\newcommand{\steleven}{\left[ {m+1} \atop {k} \right] }
\newcommand{\sttwelve}{\left[ {m+1} \atop {m-1} \right] }
\newcommand{\stthirteen}{\left[ {m+1} \atop {1} \right] }
\newcommand{\nodiv}{{|\kern-3.5pt/}}

\def\Tilde{\char126\relax}

\newtheorem{Definition}{\bf Definition}[section]
\newtheorem{Thm}[Definition]{\bf Theorem}
\newtheorem{Lem}[Definition]{\bf Lemma}
\newtheorem{Con}[Definition]{\bf Conjecture}
\newtheorem{Cor}[Definition]{\bf Corollary}

\newtheorem{Pro}[Definition]{\bf Problem}
\numberwithin{equation}{section}

\begin{document}

\title[The $2$-adic valuation of the
coefficients of a polynomial]{The $2$-adic valuation of the
coefficients of a polynomial}

\author{George Boros, V\'{\i}ctor Moll and Jeffrey Shallit}

\address{George Boros \\
Department of Mathematics, University of New Orleans \\ New
Orleans, LA 70148 \\ email: gboros@@math.uno.edu}

\address{Victor H. Moll\footnote{webpage: www.math.tulane.edu:80/\Tilde vhm} \\
Department of Mathematics, Tulane University \\ New Orleans, LA
70118 \\ email: gboros@@math.uno.edu}

\address{Jeffrey Shallit \\
Department of Computer Science, University of Waterloo \\
Waterloo, Ontario N2L 3G1  \\ Canada \\
email:shallit@@graceland.uwaterloo.ca}



\maketitle

\begin{abstract}
 In this paper we compute the $2$-adic valuations of some
polynomials associated with the definite integral $$ \int_0^\infty
{{dx} \over {(x^4+2ax^2+1)^{m+1}}} .$$ \vskip 10pt
\end{abstract}

\section{Introduction} \label{S:intro}

In this paper we present a study of the coefficients of a polynomial
defined in terms of the definite integral
\ba
N_{0,4}(a;m) & = & \int_{0}^{\infty}
\frac{dx}{(x^{4} + 2ax^{2} + 1)^{m+1}}
\label{int2}
\ea
where $m$ is a positive integer and $a > -1$ is a real number.

Apart from their intrinsic interest, these polynomials form the basis of a
new algorithm for the definite integration of rational functions.

An elementary calculation shows that
\ba
P_{m}(a) & := & \frac{2^{m+3/2}}{\pi} (a+1)^{m+1/2} N_{0,4}(a;m) \label{poly1}
\ea
\no
is a polynomial of degree $m$ in $a$ with rational coefficients.
Let
\ba
P_m (a) = \sum_{l=0}^m d_l (m) a^l.  \label{poly44}
\ea
\no
Then it can be shown that
$d_{l}(m)$ is equal to
\ba
\sum_{j=0}^{l} \sum_{s = 0 }^{m - l}
\sum_{k= s + l}^{m}
 (-1)^{k - l -s } 2^{-3k} \binom{2k}{k} \binom{2m+1}{2(s+j)}
 \binom{m- s - j}{m -k}
\binom{s+j}{j} \binom{k - s - j}{l -j } \nn
\ea

\no
from which it follows that $d_{l}(m)$ is a {\em rational number} with only a
power of $2$ in its denominator. Extensive calculations have shown that, with
rare exceptions, the numerators of $d_{l}(m)$ contain a single large prime
divisor and its remaining factors are very small. For example
\ba
d_{6}(30) & = & 2^{12} \cdot 7 \cdot 11 \cdot 13 \cdot 17 \cdot 31
\cdot 37 \cdot
639324594880985776531. \nn
\ea
\no
Similarly, $d_{10}(200)$ has $197$ digits with a prime factor
of length $137$ and its
second largest divisor is $797$. This observation lead us to investigate the
arithmetic properties of $d_l (m)$.  In this paper we discuss the $2$-adic
valuation of these $d_{l}(m)$.

The fact that the coefficients of $P_m (a)$
are {\em positive} is less elementary. This
follows from a hypergeometric representation of $N_{0,4}(a;m)$ that
implies the expression
\ba
d_{l}(m) & = & 2^{-2m} \sum_{k=l}^{m} 2^{k} \binom{2m-2k}{m-k} \binom{m+k}{m}
\binom{k}{l}. \label{simpld0}
\ea
\no
We have produced a proof of (\ref{simpld0}) that is independent of this
hypergeometric connection and is based on the Taylor expansion
\ba
\sqrt{a + \sqrt{1+c}} & = & \sqrt{a+1}
\left( 1 + \sum_{k=1}^{\infty} \frac{(-1)^{k-1}}{k} \,
\frac{P_{k-1}(a) }{2^{k+1} \; (a+1)^{k} } \; c^{k} \right); \label{taylor}
\ea
\no
see \cite{hyper} for details.

The expression (\ref{simpld0}) can be used to efficiently compute the
coefficients $d_{l}(m)$ when $l$ is large relative to $m$.
In Section~\ref{general-sec}
we derive a representation of the form
\ba
d_{l}(m) & = & \frac{1}{l!m! 2^{m+l}}
\left( \alpha_{l}(m) \prod_{k=1}^{m} (4k-1) -
\beta_{l}(m) \prod_{k=1}^{m}(4k+1) \right) \nn
\ea
\no
where $\alpha_{l}(m)$ and $\beta_{l}(m)$ are polynomials in $m$ of degrees
$l$ and $l-1$ respectively. For example
\ba
d_{1}(m) & = & \frac{1}{m! 2^{m+1}}
\left( (2m+1) \prod_{k=1}^{m} (4k-1) -
 \prod_{k=1}^{m}(4k+1) \right). \label{d1ofm}
\ea
\no
This representation can now be used to efficiently examine the coefficients
$d_{l}(m)$ when $l$ is small compared to $m$. In Section~\ref{linear-sec}
we prove that
\ba
\nu_{2}(d_{1}(m)) & = &  1-2m +  \nu_{2} \left( \binom{m+1}{2}
\right) + s_{2}(m)  \nn
\ea
\no
where $s_{2}(m)$ is the sum of the binary digits of $m$.

\section{The polynomial $P_{m}(a)$.}

Let
\ba
N_{0,4}(a;m) & = & \int_{0}^{\infty} \frac{dx}
{\left( x^{4} + 2ax^2 + 1 \right)^{m+1}} .
\nn
\ea
Then
\ba
P_{m}(a) & = & \frac{2^{m+3/2}}{\pi} (a+1)^{m+1/2} N_{0,4}(a;m) \label{poly11}
\ea
\no
is a polynomial in $a$ with positive rational coefficients.  The proof is
elementary and is presented in \cite{hyper}. It is based on the change of
variables $x = \tan \theta$ and $u = 2 \theta$ that yields
\ba
N_{0,4}(a;m) & = & 2^{-m-1} \int_{0}^{\pi}
\frac{ (1 + \cos u)^{2m+1} }{ \left( (1+a) + (1-a) \cos^{2}u \right)^{m+1}}
\; du. \nn
\ea
\no
Expanding the numerator and employing the standard substitution $z = \tan u$
produces
$$
N_{0,4}(a;m) =
$$
\begin{equation}
2^{-2m-3/2} \sum_{\nu=0}^{m} \binom{2m+1}{2 \nu}
\frac{(a-1)^{m - \nu}} {(a+1)^{m - \nu +1/2} }
\sum_{k=0}^{m- \nu} \binom{m- \nu}{k}
\frac{2^{k}}{(a-1)^{k}} B(m-k+1/2,1/2)
\label{expan3}
\end{equation}
where $B$ is Euler's beta function, defined by
$$
B(x, y) = {{\Gamma(x)\Gamma(y)} \over {\Gamma(x+y)}}.
$$
\no
The expression (\ref{poly11}) now produces the first formula for $d_{l}(m)$
given in the Introduction.

\section{The triple sum for $d_{l}(m)$.}

The expression for the coefficients $d_{l}(m)$ given in the Introduction can be
written as

\ba
\sum_{j=0}^{l} \sum_{s = 0 }^{m - l}
\sum_{k= s + l}^{m}
 (-1)^{k - l -s } 2^{-3k} \binom{2k}{k} \binom{2m+1}{2(s+j)}
 \binom{m- s - j}{m -k}
\binom{s+j}{j} \binom{k - s - j}{l -j }. \nn  \\
    & & \label{dlm}
\ea
\no
This expression follows directly from expanding (\ref{expan3})  and the value
\ba
B(j+1/2,1/2) & = & \frac{\pi}{2^{2j}} \binom{2j}{j}. \nn
\ea

\no
It follows that $d_{l}(m)$ is a rational number whose denominator is a
power of $2$, therefore
\no
\begin{Lem}
Let $p$ be an odd prime. Then
\ba
\nu_{p}(d_{l}(m)) & \geq & 0. \nn
\ea
\end{Lem}

\no
The positivity of $d_{l}(m)$ remains to be seen. \\

\section{The single sum expression for $d_{l}(m)$.}

An alternative form of the coefficients $d_{l}(m)$ is obtained by
recognizing $N_{0,4}(a;m)$ as a hypergeometric integral. A standard argument
shows that
\ba
N_{0,4}(a;m) & = & \frac{\pi \binom{2m}{m} } {2^{m+3/2} (a+1)^{m+1/2} }
{_{2}F_{1} \left[-m,m+1;1/2-m;(1+a)/2 \right]} \nn
\ea
where $_{2}F_{1}$ is a hypergeometric function, defined by
\ba
{_{2}F_{1}}[a,b,c;z] & := & \sum_{k=0}^{\infty} \frac{(a)_{k} (b)_{k} }
{ (c)_{k} k!} z^{k}, \nn
\ea
\no
where  $(r)_{k}$ is the rising factorial
\ba
(r)_{k} & = & r(r+1)(r+2) \cdots (r+k-1). \nn
\ea

\no
It follows that $P_{m}(a)$ is the {\em Jacobi polynomial} of degree $m$
with {\em parameters} $m+1/2$ and $-(m+1/2)$. Therefore
the coefficients are given by

\ba
d_{l}(m) & = & 2^{-2m} \sum_{k=l}^{m} 2^{k} \binom{2m-2k}{m-k} \binom{m+k}{m}
\binom{k}{l} \label{simpld1}
\ea
\no
from which their positivity is obvious. We have obtained a proof of
(\ref{simpld1}) that is independent of hypergeometric considerations and
is based on the presence of $P_{m}(a)$ in the Taylor expansion (\ref{taylor}).
See \cite{hyper} for details. \\

The formula (\ref{simpld1}) is very efficient for the calculation of
the coefficients $d_{l}(m)$ when $l$ approximately equal to $m$.  For
instance, we have
\ba
d_{m}(m) & = & 2^{-m} \binom{2m}{m} ; \nn \\
d_{m-1}(m) & = & 2^{-(m+1)} \binom{2m}{m}. \nn
\ea

\no
The expression (\ref{simpld1}), rewritten in the form
\ba
d_{l}(m) & = & 2^{-(2m-l)} \sum_{k=l}^{m} 2^{k-l} \binom{2m-2k}{m-k}
\binom{m+k}{m}
\binom{k}{l}, \nn
\ea
shows that
\ba
\nu_{2}(d_{l}(m)) & \geq & l-2m. \label{valuel}
\ea

\section{Basics on valuations}

Here we describe what is required on valuations. \\

\no
Given a prime $p$ and a rational number $r$, there exist
unique integers $a, b, m$
with $p \nodiv a, b$ such that
\ba
r & = & \frac{a}{b} p^{m}.
\ea
\no
The integer $m$ is the $p-$adic valuation of  $r$ and we denote it by
$\nu_{p}(r)$.  \\

\no
Now recall a basic result of number theory which states that
\ba
\nu_{p}(m!) & =
& \sum_{k=1}^{\infty} \left\lfloor \frac{m}{p^{k}} \right\rfloor.
\label{fact}
\ea
\no
Naturally the sum is finite and we can
end it at $k = \lfloor \log_{p} m \rfloor$.

\no
There is a famous result of Legendre \cite{gkp,leg}
for the $p-$adic valuation of $m!$. It
states that
\ba
\nu_{p}(m!) & = & \frac{m - s_{p}(m)}{p-1} \label{legen}
\ea
\no
where $s_{p}(m)$ is the sum  of the base$-p$ digits of $m$. In particular
\ba
\nu_{2}(m!) & = & m - s_{2}(m). \label{legen1}
\ea

\section{The constant term.}

The calculation of the $2$-adic valuation of the coefficients can be made
very explicit for the first few. We begin with the case of the constant term.

\no
We first compute
\ba
N_{0,4}(0;m) & = & \int_{0}^{\infty} \frac{dx}{(x^{4} + 1)^{m+1}} \nn
\ea
\no
via the change of variable $u=x^4$, yielding
\ba
N_{0,4}(0;m) & = & \frac{1}{4}B(1/4,m+3/4) \nn \\
& = &  \frac{\pi}{m!2^{2m+3/2}}
\prod_{k=1}^{m}(4k-1). \nn
\ea
\no
Therefore
\ba
d_{0}(m) & = & \frac{1}{m!  2^{m} } \prod_{k=1}^{m} (4k-1). \label{dzerom}
\ea
\no
\begin{Thm}
The $2$-adic valuation of the constant term $d_{0}(m)$ is given by
\ba
\nu_{2}(d_{0}(m)) & = & - \left( m + \nu_{2}(m!) \right) \nn  \\
                  & = & s_{2}(m) - 2m . \nn
\ea
\end{Thm}

\no
\begin{proof}
This follows directly from (\ref{dzerom}). The second expression comes from
(\ref{legen1}).
\end{proof}

\no
Using the single sum formula for $d_{0}(m)$ we obtain

\no
\begin{Cor}
\ba
\nu_{2} \left( \sum_{k=0}^{m} 2^{k} \binom{2m-2k}{m-k} \binom{m+k}{m} \right)
& = & m - \nu_{2}(m!) \nn \\
& = & s_{2}(m). \nn
\ea
\no
\end{Cor}

\no
\begin{Cor}
The $2$-adic valuation of the constant term $d_0 (m)$ satisfies
\ba
\nu_{2}(d_{0}(m)) & \geq & 1-2m \nn
\ea
\no
with equality if and only if $m$ is a power of  $2$.
\end{Cor}

\no
We now present a different proof of Corollary $3$ that is based on the
expression
\ba
d_{0}(m) & = & \frac{1}{m!  2^{m} } \prod_{k=1}^{m} (4k-1) \label{dzero}
\ea
\no
and the single sum formula
\ba
2^{2m} d_{0}(m) & = & \sum_{k=0}^{m} 2^{k} \binom{2m-2k}{m-k} \binom{m+k}{m}
\nn \\
& = & \binom{2m}{m} + 2
\sum_{k=1}^{m} 2^{k-1} \binom{2m-2k}{m-k} \binom{m+k}{m}. \label{single0}
\ea

\begin{proof}
     From (\ref{single0}) it follows that
\ba
\nu_{2}(d_{0}(m)) & \geq & 1-2m \nn
\ea
\no
because the central binomial coefficient is an even number.
Now from (\ref{dzero}) we obtain
\ba
\nu_{2}(d_{0}(m)) & = & -(m + \nu_{2}(m!)).
\label{ook}
\ea
\no
     From (\ref{fact}) we have
\ba
\nu_{2}(m!) & = & \sum_{k=1}^{\infty} \left\lfloor
\frac{m}{2^{k}} \right\rfloor. \nn
\ea
\no
Thus, from (\ref{ook}),
\ba
\nu_{2}(d_{0}(m)) & = & -\sum_{k=0}^{\infty} \left\lfloor
\frac{m}{2^{k}} \right\rfloor. \nn
\ea
We know $\nu_{2}(d_{0}(m)) \geq 1-2m$, so it suffices to determine
when equality occurs.  Indeed, the equation
\ba
\sum_{k=0}^{\infty} \left\lfloor \frac{m}{2^{k}} \right\rfloor
& = & 2m-1 \label{factorial}
\ea
\no
can be solved explicitly. Write $m = 2^{e}r$ with $r$ odd, and say
$2^{N} < r < 2^{N+1}$. Then
\ba
\sum_{k=0}^{\infty} \left\lfloor \frac{m}{2^{k}} \right\rfloor & = &
2^{e} \cdot r + 2^{e-1} \cdot r + \cdots + r +
\left\lfloor \frac{r}{2} \right\rfloor +
\left\lfloor \frac{r}{2^{2}} \right\rfloor +  \cdots
\left\lfloor \frac{r}{2^{N}} \right\rfloor \nn
\ea
\no
and (\ref{factorial}) leads to
\ba
r-1 & = & \sum_{k=1}^{N} \left\lfloor \frac{r}{2^{j}} \right\rfloor <
\sum_{k=1}^{N} \frac{r}{2^{j}} \leq
\sum_{k=1}^{\infty} \frac{r}{2^{j}} = \frac{r}{2} \nn
\ea
\no
and we conclude that $r=1$. The proof is finished.
\end{proof}

\section{The linear term. }
\label{linear-sec}

     From the triple sum we obtain
\ba
d_{1}(m) & = & \sum_{s=0}^{m-1} \sum_{k=s+1}^{m} (-1)^{k-s-1} 2^{-3k} (m-s)
\binom{2k}{k} \binom{2m+2}{2s+1} \binom{m-s-1}{m-k}. \nn
\ea
\no
Differentiating (\ref{poly11}) and $d_{1}(m) = P_{m}'(0)$ we produce
\ba
d_{1}(m) & = & \frac{1}{m!2^{m+1}}
\left( (2m+1) \prod_{k=1}^{m} (4k-1) - \prod_{k=1}^{m} (4k+1) \right). \nn
\ea

Therefore the linear coefficient is given in terms of
\ba
A_{1}(m) & := & (2m + 1) \prod_{k=1}^{m} (4k-1) -
\prod_{k=1}^{m}(4k+1) \label{a1ofm}
\ea
\no
so that
\ba
d_{1}(m) & = & \frac{A_{1}(m)}{m! 2^{m+1}}. \label{d1ofmnew}
\ea
We prove

\begin{Thm}
The $2$-adic valuation of the linear coefficient $d_{1}(m)$ is given
by
\ba
\nu_{2}(d_{1}(m)) & = &  1-2m +  \nu_{2} \left( \binom{m+1}{2}
\right) + s_{2}(m) . \nn
\ea
\end{Thm}

\no
Recall that the inequality $\nu_{2}(d_{1}(m)) \geq 1-2m$ follows directly
from the single sum expression.  The
theorem determines the exact value of the correction term.

\begin{proof}
We prove
\ba
\nu_{2}\left(A_{1}(m) \right) & = & \nu_{2}(2m(m+1)) \nn \\
         & = & 2 + \nu_{2} \left( \binom{m+1}{2} \right). \nn
\ea
\no
The result then follows from (\ref{legen1}) and (\ref{d1ofmnew}).

Define
\ba
B_{m} & = & \prod_{k=1}^{m} (4k+1) - 1  \nn
\ea
\no
and
\ba
C_{m} & = & (2m+1) \prod_{k=1}^{m} (4k-1) - 1. \nn
\ea
\no
Then evidently $A_{1}(m) = B_{m} - C_{m}$.

\no
We show
\no
\ba
a) \qquad \nu_{2}(B_{m}) & = & 2 + \nu_{2} \left( \binom{m+1}{2} \right) \nn \\
b) \qquad \nu_{2}(C_{m}) & \geq & 3 + \nu_{2} \left( \binom{m+1}{2} \right) \nn
\ea
\no
from which the result follows immediately.

\no
a) We have
\ba
B_{m} & = & \prod_{k=1}^{m} (4k+1) - 1  \nn  \\
      & = & \left( \sum_{j=1}^{m+1} 4^{m+1-j} \stone \right) - 1 \nn \\
      & = & \sum_{j=1}^{m} 2^{2(n+1-j)} \stone \nn \\
      & = & 2^{2} \sttwo + \sum_{k=2}^{m} 2^{2k} \stthree \nn \\
      & = & 2^{2} \binom{m+1}{2} + \sum_{k=2}^{m} 2^{2k} \stthree \nn
\ea
\no
where $\stfour$ is an (unsigned) Stirling numbers of the first kind, i.e.,
\ba
x(x+1) \cdots (x+m-1) & = & \sum_{k=0}^{m} \stfour x^{k}. \nn
\ea
\no
To prove $a)$, it suffices to show that
\ba
\nu_{2}\left( 2^{2} \binom{m+1}{2} \right) & < &
\nu_{2} \left( 2^{2k} \stthree \right) \nn
\ea
\no
for $2 \leq k \leq m$.

\no
To do this we observe that there exist integers $C_{k,i} \; (k\geq 1, \;
i \geq 0)$ such that
\ba
\stfive & = & \sum_{i=0}^{k-1} \binom{m}{2k-i} C_{k,i} \nn
\ea
see \cite[p.\ 152]{jor}. For example
\ba
\stsix & = & \binom{m}{2} \nn \\
\stseven & = & 3 \binom{m}{4} + 2 \binom{m}{3}  \nn \\
\steight & = & 15 \binom{m}{6} + 20 \binom{m}{5} + 6 \binom{m}{4} \nn \\
\stnine & = & 105 \binom{m}{8} + 210 \binom{m}{7} + 130 \binom{m}{6} +
24 \binom{m}{5}. \nn
\ea
\no
Hence the rational number
\ba
u := \frac{m(m-1) \cdots (m-k)}{(2k)!} \nn
\ea
\no
divides $\stfive$ in the sense that the quotient
\ba
{\stfive} \over u \nn
\ea
\no
is an integer.

It follows that
\ba
\nu_{2} \left( \stfive \right) & \geq & \nu_{2}( m(m-1) \cdots (m-k) )
- \nu_{2}( (2k)!) \nn \\
 & = & \nu_{2}( m(m-1) \cdots (m-k) ) -2k + s_{2}(k) \nn
\ea
\no
where  we have used (\ref{legen}).

\no
Hence, provided $k \geq 3$,
\ba
\nu_{2} \left(\stthree \right) & \geq
& \nu_{2}( (m+1)m(m-1) \cdots (m+1-k) ) - 2k + s_{2}(k)
\nn
\ea
\no
so that
\ba
\nu_{2} \left(2^{2k} \stthree \right) & \geq & \nu_{2}( (m+1)m) +
       \nu_{2}((m-1)(m-2)) + s_{2}(k)  \nn \\
   & \geq & \nu_{2}((m+1)m) + 1 + 1 \nn \\
    & > & \nu_{2} \left( 2^{2} \binom{m+1}{2} \right) \nn
\ea
\no
provided $m \geq 3$. ( For $m=1,\; 2$ it is easy to check
$\nu_{2}(B_{m}) =2$.)

On the other hand, if $k=2$, then
\ba
\stseven & = &  3 \binom{m}{4} + 2 \binom{m}{3} \nn \\
     & = & \frac{1}{24} m(m-1)(m-2)(3m-1), \nn
\ea
\no
so if $m$ is even, $m \geq 4$, we have
\ba
\nu_{2} \left( \stseven \right) & = & \nu_{2} \left( \frac{m(m-1)}{2} \right) +
   \nu_{2}(m-2) - \nu_{2}(12) \nn \\
   & \geq & \nu_{2} \left( \frac{m(m-1)}{2} \right) + 1 - 2  \nn \\
   & = & \nu_{2} \left( \frac{m(m-1)}{2} \right) - 1  \nn
\ea
\no
while if $m$ is odd, $m \geq 3$, we have
\ba
\nu_{2} \left( \stseven \right) & = & \nu_{2} \left( \frac{m(m-1)}{2} \right) +
   \nu_{2}(3m-1) - \nu_{2}(12) \nn \\
   & \geq & \nu_{2} \left( \frac{m(m-1)}{2} \right) + 1 - 2  \nn \\
   & = & \nu_{2} \left( \frac{m(m-1)}{2} \right) - 1  \nn
\ea
\no
so in either event
\ba
\nu_{2} \left( \stten \right) & \geq & \nu_{2} \left( \binom{m+1}{2} \right) -1.
\nn
\ea

Hence
\ba
\nu_{2} \left( 2^{4} \stten \right) & \geq & \nu_{2} \left( \binom{m+1}{2}
\right) + 3 \nn \\
 & > & \nu_{2} \left( 2^{m} \binom{m+1}{2} \right) \nn
\ea
\no
as desired. \\

\no
We now prove b):
\ba
C_{m} & = & (2m+1) \prod_{k=1}^{m} (4k-1) - 1. \nn
\ea
\no
We have
\ba
\prod_{k=1}^{m} (4k-1) & = & 4^{m} \prod_{k=1}^{m} (k - 1/4) \nn \\
  & = & -4^{m+1} \sum_{k=0}^{m+1} \steleven (-1/4)^{k} \nn \\
  & = & (-1)^{m} \sum_{k=1}^{m+1} \steleven (-4)^{m+1-k} \nn \\
  & = & (-1)^{m} \sum_{k=1}^{m+1} \steleven (-4)^{m+1-k} \nn
\ea
\no
thus
\ba
C_{m} &=&
\left( (-1)^{m} (2m+1) \sum_{k=1}^{m+1} \steleven
(-4)^{m+1-k}
   \right) -1. \nn
\ea
\no
When $m$ is even, we have
\ba
C_{m} & = &  \nn
\ea
\ba (2m+1) - (2m+1) \cdot 4 \sttwo -1  +
  (2m+1) \sum_{k=2}^{m} \stthree (-4)^{k} \nn
\ea
\ba
  & = & -2m^{2}(2m+3) +
  (2m+1) \sum_{k=2}^{m} \stthree (-4)^{k} \nn
\ea
\no
so, as in the proof of a), we have
\ba
\nu_{2} \left(C_{m} \right) & \geq &  \nn
\ea
\ba
\text{min}
\left( \nu_{2}(2m^{2}) ,
  \nu_{2} \left( 4^{2} \stten \right),
  \nu_{2} \left( 4^{4} \sttwelve \right),  \cdots,
  \nu_{2} \left( 4^{m} \stthirteen \right) \; \right) \nn
\ea
\ba
 & \geq & \text{min} \left( 1 + 2 \nu_{2}(m), \; 3 + \nu_{2} \left(
   \binom{m+1}{2} \right) \right) \nn
\ea
\ba
 & \geq & 3 + \nu_{2} \left( \binom{m+1}{2} \right) \nn
\ea
\no
since $m$ is even.

On the other hand, when $m$ is odd we observe that
\ba
C_{m} + 1 & = & (2m+1) \prod_{k=1}^{m} (4k-1) \nn
\ea
\no
and
\ba
C_{m+1} + 1 & = & (2m+3)(4m+3) \prod_{k=1}^{m} (4k-1) \nn
\ea
\no
so
\ba
\frac{C_{m+1} + 1}{(2m+3)(4m+3)} & = & \frac{C_{m} + 1}{2m+1} \nn
\ea
\no
and hence
\ba
C_{m} & = & \frac{(C_{m+1} + 1)(2m +1)}{(2m+3)(2m+3)} - 1 \nn \\
   & = & \frac{(2m+1)C_{m+1} - 8(m+1)^{2} }{(2m+3)(4m+3)}
\ea
\no
so
\ba
\nu_{2}(C_{m}) & \geq & \text{min} \left( \nu_{2}(C_{m+1}),
   2 \nu_{2}(m+1) + 3 \right) \nn \\
  & \geq & 3 + \nu_{2} \left( \binom{m+1}{2} \right) \nn
\ea
\no
since $m$ is odd.

\no
This completes the proof.
\end{proof}

\no
The corresponding question of the $3$-adic valuation of $d_{1}(m)$ seems to be
more difficult. We  propose

\begin{Pro}
Prove the existence of a sequence of positive integers $m_{j}$ such that
$\nu_{3}(d_{1}(m_{j})) = 0$.  Extensive calculations show that
\ba
m_{j+1} - m_{j}  & \in & \left\{ 2, \, 7, \, 20, \, 61, \, 182, \cdots \right\}
\label{16sequence}
\ea
\no
where the sequence $ \{ q_{j} \}$ in (\ref{16sequence})
is defined by $q_{1} = 2$ and
$q_{j+1} = 3q_{j} + (-1)^{j+1}$. It would be of interest to know whether
$\nu_{3}(d_{1}(m))$ is unbounded: the maximum value for $2 \leq m \leq
20000$ is $12$, so perhaps $\nu_{3}(d_{1}(m)) = O( \log m)$ as
$m \to \infty$.
\end{Pro}

\section{The general situation.}
\label{general-sec}

In this section we prove the existence of  polynomials
$\alpha_{l}(x)$ and $\beta_{l}(x)$
with positive integer coefficients such that
\ba
d_{l}(m) & = & \frac{1}{l!m! 2^{m+l}}
\left( \alpha_{l}(m) \prod_{k=1}^{m} (4k-1) -
\beta_{l}(m) \prod_{k=1}^{m}(4k+1) \right). \nn
\ea

\no
These  polynomials are efficient for the calculation of $d_{l}(m)$ if $l$ is
small relative to $m$, so they  complement the results of Section $4$.

\no
For example
\ba
\alpha_{0}(m) & = & 1 \nn \\
\alpha_{1}(m) & = & 2m+1 \nn \\
\alpha_{2}(m) & = & 2(2m^2 + 2m+1) \nn \\
\alpha_{3}(m) & = & 4(2m+1)(m^2 + m + 3) \nn \\
\alpha_{4}(m) & = & 8(2m^{4} + 4m^{3} + 26m^{2} + 24m + 9). \nn
\ea
\no
and
\ba
\beta_{0}(m) & = & 0 \nn \\
\beta_{1}(m) & = & 1 \nn \\
\beta_{2}(m) & = & 2(2m+1) \nn \\
\beta_{3}(m) & = & 12(m^2 + m + 1) \nn \\
\beta_{4}(m) & = & 8(2m+1)(2m^2+2m+9). \nn
\ea

\no
The proof consists in computing the expansion of $P_{m}(a)$ via the
Leibnitz rule:
\ba
P_{m}(a) & = & \frac{2^{m+3/2}}{\pi}
\sum_{j=0}^{l} \binom{l}{j} \left( \frac{d}{da} \right)^{l-j} (a+1)^{m+1/2}
\Big|_{a=0}
\left( \frac{d}{da} \right)^{j} N_{0,4}(a;m) \Big|_{a=0}. \nn
\ea
\no
We have
\ba
\left( \frac{d}{da} \right)^{r} (a+1)^{m+1/2}
\Big|_{a=0} & = & 2^{-2r} \frac{(2m+2)!}{(m+1)!}
\frac{(m-r+1)!}{(2m-2r+2)!} \label{derivpower}
\ea
\no
and
\ba
\left( \frac{d}{da} \right)^{r} N_{0,4}(a;m) \Big|_{a=0} & = &
(-1)^{r} \frac{(m+r)!}{m!} 2^{r} \int_{0}^{\infty}
\frac{x^{2r} }{(x^{4} + 1)^{m+r+1}} \; dx. \label{derivinteg}
\ea
\no
The integral is evaluated via the change of variable $t=x^4$ as
\ba
\int_{0}^{\infty}
\frac{x^{2r} \; dx }{(x^{4} + 1)^{m+r+1}} & = &
\frac{1}{4} B \left( \frac{r}{2} + \frac{1}{4}, m + \frac{r}{2} + \frac{3}{4}
\right). \nn
\ea
\no
This yields
\ba
\left( \frac{d}{da} \right)^{r} N_{0,4}(a;m) \Big|_{a=0} & = &
\frac{(-1)^{r} (2r)!}{2^{2r+2m+3/2}} \; \frac{\pi}{m!r!}
\prod_{l=1}^{m}(4l-1+2r). \label{final}
\ea
\no
Therefore
\ba
P_{m}^{(l)}(0) & = & \nn
\ea
\ba
\frac{l! (2m+2)!}{2^{m+2l}m!(m+1)!}
\sum_{j=0}^{l}
\frac{(-1)^{j} (m-l+j+1)! (2j)!}
{j!^{2} (l-j)!(2m-2l+2j+2)!}
\prod_{\nu=1}^{m} (4 \nu -1 + 2j). \nn
\ea
\no
We now split the sum according to the parity of $j$. In the case
$j$ is odd $( = 2t-1 )$ we use
\ba
\prod_{\nu=1}^{m}(4 \nu -1 + 2j) & = &
\prod_{\nu=1}^{m}(4 \nu +1)
\left( \prod_{\nu=m+1}^{m+t-1}(4 \nu +1) \Big{/}
\prod_{\nu=1}^{t-1}(4 \nu +1) \right) \nn
\ea
\no
and if $j$ is even ($ = 2t$) we employ
\ba
\prod_{\nu=1}^{m}(4 \nu -1 + 2j) & = &
\prod_{\nu=1}^{m}(4 \nu -1)
\left( \prod_{\nu=m+1}^{m+t}(4 \nu -1) \Big{/}
\prod_{\nu=1}^{t}(4 \nu -1) \right). \nn
\ea
\no
We conclude that
\ba
d_{l}(m) & = & X(m,l) \prod_{\nu=1}^{m} ( 4 \nu -1) - Y(m,l)
\prod_{\nu=1}^{m} (4 \nu +1) \nn
\ea
\no
with
\ba
X(m,l) & = & \nn
\ea
\ba
\frac{(2m+2)!}{2^{m+2l} m! (m+1)!}
\sum_{t=0}^{\lfloor l/2 \rfloor}
\frac{(m-l+2t+1)! (4t)!}{(2t)!^{2} (l-2t)! (2m-2l+4t+2)!}
{{\prod_{\nu=m+1}^{m+t} (4 \nu -1)} \over {\prod_{\nu=1}^{t} (4 \nu -1)}} \nn
\ea
\no
and
\ba
Y(m,l) & = & \nn
\ea
\ba
 \frac{(2m+2)!}{2^{m+2l} m! (m+1)!}
\sum_{t=1}^{\lfloor (l+1)/2 \rfloor}
\frac{(m-l+2t)! (4t-2)!}{(2t-1)!^{2} (l-2t+1)! (2m-2l+4t)!}
{{ \prod_{\nu=m+1}^{m+t-1} (4 \nu +1)} \over
{\prod_{\nu=1}^{t-1} (4 \nu +1)}}. \nn
\ea
\no
The quotients of factorials appearing above can be simplified via
\ba
\frac{(m+1)!}{(m-l+2t+1)!} & = & \prod_{j=1}^{l-2t} (j+m-l+2t+1) \nn
\ea
\no
and
\ba
\frac{(2m+2)!}{(2m-2l+4t+2)} & = & 2^{l-2t}
\left(\prod_{i=1}^{l-2t} (i+m-l+2t+1) \right)
\left( \prod_{i=1}^{l-2t} (2i+2m-2l+4t+1) \right). \nn
\ea
\no
We conclude that
\ba
d_{l}(m) & = & \frac{1}{l!m!2^{m+l}}
\left( \alpha_{l}(m) \prod_{\nu=1}^{m} ( 4 \nu -1) -
\beta_{l}(m) \prod_{\nu=1}^{m} (4 \nu +1) \right) \nn
\ea
\no
with
\ba
\alpha_{l}(m) & = & l! \sum_{t=0}^{\lfloor l/2 \rfloor}
\frac{\binom{4t}{2t} }{2^{2t} (l-2t)! }
\frac{\prod_{\nu=m+1}^{m+t}}{\prod_{\nu=1}^t (4 \nu -1) }
\left( \prod_{\nu=1}^{t} (4 \nu -1) \right)
\left( \prod_{\nu=m-(l-2t-1)}^{m} ( 2 \nu +1) \right) \nn
\ea
\no
and
\ba
\beta_{l}(m) & = & l! \sum_{t=1}^{\lfloor (l+1)/2 \rfloor}
\frac{\binom{4t-2}{2t-1} }{2^{2t-1} (l-2t+1)! }
\left( \frac{\prod_{\nu=m+1}^{m+t-1}
(4 \nu +1) }{\prod_{\nu=1}^{t-1} (4 \nu +1) } \right)
\left( \prod_{\nu=m-(l-2t)}^{m} ( 2 \nu +1) \right). \nn
\ea
\no
The identity
\ba
\prod_{\nu=1}^{t} (4 \nu -1) & = & \frac{(4t)!}{2^{2t}(2t)!}
\; \left( \prod_{\nu=1}^{t-1} (4 \nu +1) \right)^{-1} \nn
\ea
is now employed to produce
\ba
\alpha_{l}(m) & = & \sum_{t=0}^{\lfloor l/2 \rfloor} \binom{l}{2t}
\prod_{\nu=m+1}^{m+t} ( 4 \nu -1)
\prod_{\nu=m-(l-2t-1)}^{m} ( 2 \nu +1)
\prod_{\nu=1}^{t-1} ( 4 \nu +1)  \nn
\ea
\no
and
\ba
\beta_{l}(m) & = & \sum_{t=1}^{\lfloor (l+1)/2 \rfloor} \binom{l}{2t-1}
\prod_{\nu=m+1}^{m+t-1} ( 4 \nu +1)
\prod_{\nu=m-(l-2t)}^{m} ( 2 \nu +1)
\prod_{\nu=1}^{t-1} ( 4 \nu -1).  \nn
\ea

We have proven:

\begin{Thm}
There exist polynomials $\alpha_{l}(x)$ and $\beta_{l}(x)$ with integer
coefficients such that
\ba
d_{l}(m) & = & \frac{1}{l!m! 2^{m+l}}
\left( \alpha_{l}(m) \prod_{k=1}^{m} (4k-1) -
\beta_{l}(m) \prod_{k=1}^{m}(4k+1) \right). \nn
\ea
\end{Thm}

\no
Based on extensive numerical calculations we propose

\begin{Con}
All the roots of the polynomials $\alpha_{l}(m)$ and $\beta_{l}(m)$ lie
on the line $\text{Re}(m) = -1/2$.
\end{Con}

\end{document}